\documentclass[12pt]{article}
\usepackage{latexsym,amsthm,amsfonts,
epsfig,
psfrag,
}
\begin{document}

\noindent

\newtheorem{lemma}{Lemma}
\newtheorem{theorem}{Theorem}
\renewcommand{\proofname}{Proof}

\def \S{{\mathbb S}}
\def \X{{\mathbb X}}
\def \H{{\mathbb H}}
\def \E{{\mathbb E}}

\begin{center}
{\Large  \bf
On finite index reflection subgroups 
of discrete reflection groups.
}

\vspace{9pt}
{\large A.~Felikson, P.~Tumarkin}
\end{center}

\vspace{4pt}
\begin{center}
\parbox{10cm}{\scriptsize
{\bf Abstract.}
Let $G$ be a  discrete group 
generated by reflections in hyperbolic or Euclidean space,
and $H\subset G$ be a finite index subgroup generated by reflections.
Suppose that the fundamental chamber of $G$ is a finite volume
polytope  with $k$ facets.
In this paper, we prove that the fundamental chamber of $H$ has at least
$k$ facets.

}\\
\end{center}

\vspace{10pt}

{\bf 1.}
Let $\X$ be hyperbolic space $\H^n$, 
Euclidean space $\E^n$ or spherical space $\S^n$.
A polytope in $\X$
is called a {\it Coxeter polytope} if all its dihedral angles are
integer parts of $\pi$.
A convex polytope in $\X$ admits a {\it Coxeter decomposition}
if it is tiled by finitely many Coxeter polytopes
such that any two tiles having a common facet are symmetric with respect
to this facet.
The tiles of the decomposition are called {\it fundamental polytopes}.


In this paper we show that 
if the polytope admitting a Coxeter decomposition has exactly $k$ facets
then the fundamental polytope has at most $k$ facets (Theorem~1).
In particular,
consider a finite covolume discrete group $G$
generated by reflections in hyperbolic or Euclidean space and a finite
index reflection subgroup $H\subset G$.
If the fundamental chamber of $G$ has exactly $k$ facets, 
then the fundamental chamber of $H$ has at least $k$ facets.

\vspace{7pt}
\noindent
The authors are grateful to E.~B.~Vinberg for useful discussions and remarks.

\vspace{10pt}

{\bf 2.} 
A polytope $P\subset \X$ is called {\it acute-angled}
if its dihedral angles are less or equal to
 $\frac{\pi}{2}$.
The minimal plane containing a face of a polytope we call
an {\it extension} of this face. 
The proof of Theorem~1 is based on the following fact
proved by E.~M.~Andreev (see~\cite{Andr2}):

{\it
Let $P\subset \X$ be an acute-angled polytope.
If two faces of $P$ have no common points then the extensions of these
faces have no common points.

}

A convex polytope with fixed Coxeter decomposition will be denoted
by $P$, the fundamental polytope of the decomposition will be denoted by $F$.
A hyperplane $H$ is called a {\it mirror} of the decomposition
 if it contains a facet 
of a fundamental polytope and contains no facet of $P$.

Let $\Theta$ be a Coxeter decomposition of $P$.
A dihedral angle of $P$ formed up by  facets $\alpha$  and $\beta$
is called {\it fundamental} in the decomposition
$\Theta$ if no mirror of  $\Theta$ contains
$\alpha \cap \beta$.

\begin{lemma}\label{intersect}
Let $P\subset X$ be a convex polytope with exactly $k$ facets.
Let $\alpha$ be a hyperplane decomposing $P$ into polytopes
$P_1$  and $P_2$. If  each of $P_1$ and $P_2$ has at least $k+1$ facets then   
$\alpha$ intersects the inner part of each facet of $P$
(in particular, $\alpha$ cuts no dihedral angle).
Furthermore, if $P_1$ is a Coxeter polytope then $P$ is a Coxeter
polytope, too.

\end{lemma}

\begin{proof}
Let $\beta_1,...,\beta_k$ be the facets of $P$.
Suppose that $\alpha$ does not intersect the inner part of $\beta_1$.
Then $\beta_1$ contains no facet of either $P_1$ or $P_2$. 
If $\beta_1$ contains no facet of $P_1$ then
$P_1$ has at most $k$
facets (one facet is contained in $\alpha$ and the rest facets are 
contained in $\beta_2,...,\beta_k$).

Suppose that $P_1$ is a Coxeter polytope.
Then the extensions of two facets of $P_1$ have a common point
if and only if the facets have a common point (\cite{Andr2}).
Each facet of $P$ contains some facet of $P_1$.
Clearly, the corresponding dihedral angles of $P$ and $P_1$ coincide.
Thus,  $P$ is a Coxeter polytope.

\end{proof}

\begin{lemma}\label{simp}
Let $P\in X$ be a convex polytope and  
$\alpha$ be a hyperplane containing an inner point of each facet of $P$. 
Suppose that $\alpha$ decomposes $P$ into polytopes
$P_1$ and $P_2$, where $P_1$ is a Coxeter polytope.
Then $\alpha$ contains no vertex of $P$.

\end{lemma}

\begin{proof}

Suppose that $\alpha$ contains a vertex $A$  of $P$.
Let $s$ be a sphere centered in $A$ (if $A$ is an ideal vertex of a
hyperbolic polytope consider a horosphere $s$ centered in $A$).
Consider a section $p$ of $P$ by $s$.
Let $p_1$ be a section of $P_1$ by $s$.
Since $\alpha$ intersects the inner part of each facet of $P$,
the set of facets of the spherical (Euclidean) polytope $p$
coincides with the set of facets of $p_1$ without the facet $\alpha\cap s$.
By Lemma~\ref{intersect} $P$ is a Coxeter polytope.
Hence, $p_1$ and $p$ are Coxeter polytopes too.
Let $\Sigma(p)$ and $\Sigma(p_1)$ be the Coxeter diagrams of 
$p$ and $p_1$ and  $v$ be a node of $\Sigma(p_1)$ corresponding to
$\alpha\cap s$. Then $\Sigma(p)=\Sigma(p_1)\setminus v$.

Suppose that $A$ is not an ideal vertex.
Then $p$ and $p_1$ are spherical simplices of the same dimension.
This is impossible, since these polytopes have different numbers of facets.
Suppose that $A$ is an ideal vertex.
Since any compact Euclidean Coxeter polytope is a direct product of simplices,
$\Sigma(p_1)$ is a disjoint union of connected parabolic diagrams.
But the subdiagram $\Sigma(p_1)\setminus v=\Sigma(p)$ of $\Sigma(p_1)$
is a disjoint union of  connected parabolic diagrams too.
This is impossible by the definition of connected parabolic diagram.
The contradiction proves the lemma.

\end{proof}

\begin{lemma}\label{non-acute}
Let $P$ be either a convex finite volume polytope in  $\H^n$ or $\E^n$ 
or a convex polytope in $\S^n$ containing no pair of antipodal points 
of the sphere.
Let 
$\alpha$ be a hyperplane containing an inner point of each facet of $P$. 
Suppose $\alpha$ decomposes $P$ into polytopes
$P_1$ and $P_2$.
If $\alpha$ contains no vertex of $P$ then neither $P_1$ nor $P_2$ 
is acute-angled.

\end{lemma}

\begin{proof}

Suppose  $P_1$ is an acute-angled polytope.
Let $H$ be a face of $P_2$ different from an ideal vertex.
Let us prove that the extension of $H$ contains a face of $P_1$.
If $H$ is a facet the statement is true by the
assumption.
Let $H_t$ be a face of codimension $t$ different from an ideal vertex. 
Suppose inductively that the statement is true for the faces of
codimension smaller than $t$.
Consider $H_t$ as an intersection of a facet $H_1$ with a face $H_{t-1}$
of codimension $t-1$ (here $H_1$ and $H_{t-1}$ are the faces of $P_2$).
By the assumption the extensions of $H_1$ and $H_{t-1}$  contain some faces
of the acute-angled polytope $P_1$.
Thus, by~\cite{Andr2} the corresponding faces of $P_1$ have a common
point. Therefore, the extension of $H_t$ contain some face of $P_1$.
We have proved that the extension of any face of $P_2$
different from an ideal vertex contains a face of $P_1$.

Let $A$ be a vertex of $P$ contained in $P_2$. 
By the hypothesis of the lemma $A$ does not belong to $\alpha$. 
Let $AB_1$,...,$AB_r$ be the edges of $P$ incident to $A$. 
We have proved that each of these edges contains a face of $P_1$.
Thus, the points $B_1$,...,$B_r$ are the vertices of $P_1$. 
Hence, $A$ is a unique common vertex of $P$ and $P_2$.
Let $H$ be any facet of $P$ that is not incident to
$A$ (this facet exists, since $P$ is a finite volume polytope
and $P$ is not spherical polytope containing a pair of antipodal points). 
Then $H$ is contained in one half-space with respect to
$\alpha$. The contradiction with the hypothesis proves the lemma.

\end{proof}

\begin{theorem}\label{dec}
Let $P$ be a finite volume convex polytope in $\H^n$ or $\E^n$
with exactly $k$ facets.
If $P$ admits a Coxeter decomposition  $\Theta$ with fundamental polytope $F$,
then $F$ has at most $k$ facets.
 
\end{theorem}

\begin{proof}
Suppose the contrary.
Let $\Theta$ be a Coxeter decomposition of $P$
with fundamental polytope $F$ bounded by more than $k$ facets.
Without loss of generality we can assume that the theorem is true for 
any polytope in $\H^n$ ($\E^n$) bounded by less than $k$ facets.

Let $M$ be a set of convex polytopes such that
if $P_1\in M$ then

 1) each facet of $P_1$ is contained either 
in a mirror of $\Theta$ or in a facet of $P$;

2) $P_1$ has exactly $k$ facets.

Since a number of mirrors in a Coxeter decomposition is finite, 
$M$ is a finite nonempty set.
Let $P_{min}\in M$ be a polytope minimal with respect to inclusion.
Let $\Theta_{min}$ be a restriction of $\Theta$ to the polytope $P_{min}$.  
Then $\Theta_{min}$ is a Coxeter decomposition with the fundamental
polytope $F$. Note that $P_{min}$ has $k$ facets and $F$ has more than $k$
facets.

Let $\phi$ be a dihedral angle of $P_{min}$.
Suppose that $\phi$ is not fundamental in  $\Theta_{min}$.
Let $\beta$ be a mirror of  $\Theta_{min}$ decomposing $\phi$. 
Then $\beta$ decomposes  $P_{min}$ into two polytopes 
$P'_{min}$ and $P''_{min}$, where each of $P'_{min}$ and $P''_{min}$
has at most $k$ facets. 
If $P'_{min}$ has exactly $k$ facets, we have a contradiction with 
the minimal property of $P_{min}$ in $M$. Hence,  $P'_{min}$ has less
than $k$ facets. Consider a restriction  $\Theta'_{min}$ of
$\Theta_{min}$ to the polytope $P'_{min}$.  $\Theta'_{min}$ is a
Coxeter decomposition with the fundamental polytope $F$, where
$F$ has more than $k$ facets. 
This contradicts to the assumption that the theorem is true for 
any polytope in $\H^n$ ($\E^n$) bounded by less than $k$ facets.
Thus, any dihedral angle of $P_{min}$ is fundamental.

Let $N$ be a set of polytopes such that  
if $P_1\in N$ then $P_1$ is bounded by exactly $k$ facets of $P_{min}$
and by a mirror of $\Theta_{min}$.
Clearly, $N$ is a finite nonempty set. 
Let $P_1\in N$ be a polytope minimal with respect to inclusion.

Let $\alpha$ be a mirror decomposing $P_{min}$ into polytopes
$P_1$ and $P_2$.
If $\phi$ is a dihedral angle of $P_1$ formed up by two facets of 
$P_{min}$ then $\phi$ is fundamental as it is already proved.
Let $\phi$ be a dihedral angle of $P_1$ formed up by  $\alpha$ and 
a facet of $P_{min}$.
Suppose that  $\phi$ is not fundamental.
Let  $\beta$ be a mirror decomposing  $\phi$. 
Suppose that $\beta$ intersects an inner part of each facet of $P_{min}$.
Then $\beta$ cuts a polytope  $P'_1\subset P_1$ from $P_{min}$ 
 such that $P'_1\in N$.
This contradicts to the minimal property of $P_1$ in $N$.
Now suppose that $\phi$ does not intersect some facets of $P_{min}$.
Then $\beta$ cuts a polytope $P'_1\subset P_1$ from $P_{min}$ 
 such that $P'_1$ has at most $k$ facets.
This is impossible by the minimal property of $P_{min}$ in $M$
and the assumption that the theorem is true for 
any polytope in $\H^n$ ($\E^n$) bounded by less than $k$ facets.

Thus, any  dihedral angle of $P_1$ is fundamental
and  $P_1$ is a Coxeter polytope.
By Lemma~\ref{simp},  $\beta$ contains no vertex of $P_{min}$. 
By Lemma~\ref{non-acute}  $P_1$ is not an acute-angled polytope.
The contradiction to the fact that any Coxeter polytope is
acute-angled proves the theorem.

\end{proof}

\noindent
{\bf Remark 1.} 
Clearly, the theorem is false for $\S^n$. 
An additional assumption that $P$  contains no pair of antipodal
points makes the theorem
true but trivial. With this assumption both polytopes $P$ and $F$ 
are simplices
(see~\cite{Cox} for the classification of spherical Coxeter polytopes). 

\vspace{6pt}
\noindent
{\bf Remark 2.}
The theorem is false without the assumption that $P$ is of finite volume.  
For example, let $G$ be a group generated by the reflections with
respect to two divergent hyperplanes
$m_1$ and $m_2$ in  $\H^n$
(or two parallel hyperplanes in $\E^n$).
Then $G$ has an index 2 subgroup 
generated by the reflections with respect  to $m_1$, $m_2$ and $h$,
where $h$ is perpendicular to $m_1$ and $m_2$.


\vspace{6pt}
\noindent
{\bf Remark 3.} 
Let $P\subset \X$ be a finite volume Coxeter polytope with $k$ facets
and $G_P$ be a group generated by the reflections with respect to the
facets of $P$.
There are some examples of polytopes such that the group $G_P$ can be 
generated by less than $k$ generators. However, in these examples 
some generators are not reflections. This leads to the conjecture
that $G_P$ can not be generated by less than $k$ reflections.
If the conjecture is true, we can prove the following corollary
of Theorem~1:

\vspace{6pt}
{\it
Let $G$ be a discrete group generated by reflections in $\H^n$ or
$\E^n$, $F$ be a fundamental chamber of $G$ and $H\subset G$
be  a finite index subgroup generated by reflections.
If $F$ is a finite volume polytope with $k$ facets then $H$
can not be generated by less than $k$ reflections.
}

\vspace{20pt}

\noindent
{\it Independent University of Moscow\\
 E-mail addresses:\quad felikson@mccme.ru\\
\phantom{E-mail addresses:}\quad pasha@mccme.ru
}


\begin{thebibliography}{3}

\bibitem{Andr2}
E.~M.~Andreev,
Intersection of plane boundaries of acute-angled polyhedra.
Math. Notes 8 (1971) 761--764.

\bibitem{Cox}
H.~S.~M.~Coxeter, Finite groups generated by reflections,
and their subgroups generated by reflections,
PRC. Cambridge Philis. Soc. 30 (1934) 466--482.\\



\end{thebibliography}
\end{document}